\documentclass{article}
\usepackage{amsfonts}
\usepackage{amsthm}
\begin{document}
\author{Jeremy Berquist}
\title{The Obstacle Problem on Graphs and Other Results}
\maketitle

\textbf{Abstract.  }  Our primary motivation is existence and uniqueness for the obstacle problem on graphs.  That is, we look for unique solutions to the problem $Lu = \chi_{\{u>0\}}$, where $L$ is the Laplacian matrix associated to a graph, and $u$ is a nonnegative real-valued vector with preassigned zero coordinates and positive coordinates to be determined.  In the course of solving this problem, we make a detour into the study of Laplacian matrices themselves.  First, we present the row reduced echelon form of such matrices and determine the invertibility of proper square submatrices.  Next, we determine eigenvalues of several simple Laplacians.  In this context, we introduce a new polynomial called the generalized characteristic polynomial that allows us to compute (theoretically, if inefficiently)the usual characteristic polynomial for trees by inspection of the graph.  Finally, we give our solution to the obstacle problem on graphs and discuss other components of the obstacle problem, which we investigate in future research.

\section{Introduction}
\noindent
An obstacle problem is a free boundary problem in which the equilibrium position of an elastic membrane is sought, when the boundary is held fixed and the membrane is constrained to lie above a given obstacle.  One formulation of an obstacle-type problem is to seek nonnegative solutions $u$ of the equation $$Lu = \chi_{\{u>0\}},$$ where $L$ is the Laplacian operator, and $u$ is a twice differentiable function; in other words, the left hand side is $$\sum \frac{\partial^2 u}{\partial x_i^2},$$ the sum being over all local coordinates $x_i$.  The right hand side is the characteristic function for the subset on which $u$ is strictly positive.  One thinks of $u$ as being the distance between the elastic membrane and the obstacle over which it is stretched.  This problem has been studied in several guises, varying the conditions on the function $u$ and possibly the operator being applied.  In the present form, it has been studied by Teka and Blank [BT].

In the discrete version of this problem, $u$ is a function of the vertices of a finite graph, the graph being connected, undirected, without loops (edges emanating from and terminating at the same vertex), and with at most one edge between vertices.  That is, $u$ is an $n$-vector with nonnegative real coordinates, each coordinate corresponding to a vertex in the graph.  The Laplacian operator in this case is the Laplacian matrix (also called the adjacency matrix) associated to the graph.  The characteristic function on the right hand side has its usual meaning.  Moreover, when we consider this equation, we are considering several different scenarios.  The problem is to assign $u_i = 0$ for some vertices, and then to determine the positive values of the other $u_i$.  For every such assignment, we ask whether the nonzero values are positive and uniquely determined.

It turns out that the discrete problem, as formulated, has only the trivial solution $$u_1 = u_2 = \cdots = u_n = 0.$$  To see this, it suffices to observe that the sum of the rows in any Laplacian matrix is the zero row.  Thus it becomes necessary either to change the problem to a form suited to graphs, or to change the nature of the vector $u$.  We investigated both changes and found ``reasonable" solutions in either case.  In the course of solving the discrete obstacle problem, it was necessary to find the rank of an arbitrary Laplacian matrix.  Questions about the row reduced echelon form and eigenvalues arose naturally, even as they are not directly related to the obstacle problem.

Existence and uniqueness are just one part of the obstacle problem on graphs.  Regularity and nondegeneracy can also be considered.  See for instance [PSU] for the definitions of these conditions for standard obstacle problems in PDEs.  Regularity of solutions involves showing that solutions to an obstacle-type problem have the same topological and boundedness properties exhibited by functions appearing in the statement of the problem.  Nondegeneracy has to do with comparing behaviors of two solutions on the free boundary (the intersection of the boundary of the set where $u>0$ with the boundary of the set where $u=0$) with their behaviors away from the free boundary.  It remains to find the correct interpretation of these properties for the discrete case.  Topologically, one can introduce a metric on a graph, where open balls consist of vertices sufficiently close to a given vertex by a shortest possible path (assuming the graph is connected).  Whether solutions to the obstacle problem on graphs have these properties deserves more investigation.

The rest of this paper is a compendium of results related to Laplacian matrices themselves.  Only in future work will we deal with the other parts of the discrete obstacle problem.  Several authors have discovered methods for finding eigenvalues of Laplacian matrices and for determining their sizes.  See for instance [S] or [M].  After investigating invertibility properties and eigenvalues in the next two sections of the paper, we turn to existence and uniqueness of solutions to the obstacle problem in the final section.  Thus the reader can safely skip the section on eigenvalues.

\section{Preliminaries}
\noindent
We will consider finite graphs without loops (edges emanating from and terminating at the same vertex) that are undirected and have at most one edge between vertices.  We will frequently label the vertices $v_1, v_2, \ldots, v_n$.  When we have occasion to specify nonnegative real numbers for each vertex, we will use $u_1, u_2, \ldots, u_n$.  Given such a graph, the Laplacian matrix $L$ associated to the graph is the $n \times n$ real symmetric matrix with entries \begin{itemize}
\item{$L_{ii}= $ the number of edges emanating from $v_i$,}
\item{$L_{ij} = -1$ if $v_i$ and $v_j$ are adjacent, and }
\item{$L_{ij} = 0$ otherwise.}
\end{itemize}
When we consider the equation $Lu = \chi_{\{u>0\}},$ we mean that $u$ is a real-valued vector with nonnegative entries, whose zero entries have been specified in advance and whose positive entries are to be determined.  We consider only the case where at least one entry is zero, and where at least one entry is positive.  This is the form of the classical obstacle problem in PDEs.  As we mentioned before, the only sensible solution is $u_1 = u_2 = \cdots = u_n = 0.$  Thus, we have to reformulate the problem in order to consider nonzero solutions.  

One way of treating the discrete problem is to ignore those equations in the system $Lu = \chi_{\{u>0\}}$ where the right hand side is 0, and to determine the positive entries of $u$ from only those equations whose right hand side is 1.  In other words, we think of $u$ as a vector whose entries $u_i$ are \textit{functions} of the vertices.  Then we require that if $u_i$ is zero at $v_i$, it is zero everywhere, and if $u_i$ is positive at $v_i$, it is positive everywhere.  We let the entries $u_i$ be vectors themselves and are only interested in their values at the corresponding vertices $v_i$.  In this way, then, the equations with 0 right hand side make sense, even if the $u_i$ appearing there are thought of as positive (at their own respective vertices).

For instance, for the path graph on 3 vertices, suppose we set $u_1 = 0$ and $u_2$ and $u_3$ are positive.  The system to solve is $$\begin{array}{c}
-u_2 = 0, \\
2u_2 - u_3 = 1, \\
u_3 - u_2 = 1. \\
\end{array}$$
Of course, this system is inconsistent.  When rewriting according to the above interpretation, the system we solve is $$\begin{array}{c}
2u_2-u_3 = 1, \\
u_3 - u_2 = 1,  \\
\end{array}$$
which has the unique solution $u_2 = 2$, $u_3 = 3.$

With this version of the obstacle problem on graphs, we have the following result, which we prove in the final section:

\noindent
\linebreak
\textbf{Theorem 4.1.  }  There exists a unique positive solution $u$ to the obstacle problem $Lu = \chi_{\{u>0\}}$.  In other words, for any assignment of zero entries to the vector $u$, the other entries are determined uniquely as positive real numbers by the nonzero equations in the system $Lu = \chi_{\{u>0\}}.$

\noindent
\linebreak
This interpretation of the obstacle problem on graphs seems to produce the right solution, but it involves a slightly unnatural interpretation of the vector $u$.  Moreover, ignoring the equations with zero right hand side involves forgetting about important properties of the original system.  To address this difficulty, we have also considered a slightly modified system $$Lu + b = \chi_{\{u>0\}},$$ where $b$ is now a random vector with entries $b_i \geq 0$, whose values are uniquely determined by the assignment of zeros, as follows:  for each $i$ such that $u_i >0$, $b_i = 0$, and the other coordinates of $b$ are uniquely determined by the equations with 0 right hand side in the given system.  Adding such an error term eliminates the fundamental inconsistency of the system of equations $Lu = \chi_{\{u>0\}}.$  

This seems to be a more natural interpretation of the obstacle problem for graphs.  We prove the following result in the final section (although the proof is essentially the same as that of Theorem 4.1).

\noindent
\linebreak
\textbf{Theorem 4.2.  }  There exists a unique solution to the obstacle problem $Lu + b = \chi_{\{u>0\}}.$  In other words, for any assignment of zero entries to the vector $u$, the random entries of $b$, as well as the positive entries of $u$ are uniquely determined.

\noindent
\linebreak
The issue with our first solution for the obstacle problem is that the positive values $u_i$ are uniquely determined and positive, but that there are some auxiliary equations that the $u_i$ should solve and unfortunately do not.  We dealt with this problem by simply ignoring those equations.  In our second solution, we have introduced a random vector, effectively saying that once the positive values of $u_i$ are found, then the auxiliary equations (those with 0 right hand side in the original obstacle problem $Lu = \chi_{\{u>0\}}$) are satisfied \textit{by definition}.  Considering again the path graph on 3 vertices, where $u_1 = 0$, the system is $$\begin{array}{c}
-u_2 = -b, \\
2u_2 - u_3 = 1, \\
u_3 - u_2 = 1.  \\
\end{array}$$
There is no inconsistency here.  The values $u_2 = 2$ and $u_3=3$ determined by the last two equations as before, and then we set $b = 2$.  

We address, briefly, in the last section a more moderate approach than these.  In particular, instead of adding a random vector, we add a constant vector.  This takes care of the fundamental inconsistency of the system of equations without allowing auxiliary equations to be true ``by definition."  It turns out that solutions exist uniquely and are positive, when they exist, but that existence depends on the particular graph being studied.  For complete graphs, there is always a solution.  For path graphs and cycle graphs, the zeros have to be placed symmetrically about the graph.

Now we turn to the other parts of the paper.

As a real symmetric matrix, $L$ is diagonalizable.  In fact, it follows directly from the definition that the sum of the rows is the zero row for any Laplacian matrix $L$.  Thus, $L$ is not invertible; equivalently, 0 is an eigenvalue.  As a first result, which will be used elsewhere in the paper, we show that certain square submatrices taken from Laplacian matrices are invertible and have inverses with only nonnegative entries.  In particular, if we choose fewer than $n$ rows from $L$, and consider the square submatrix obtained by using these rows and their corresponding columns, then the inverse has only nonnegative entries.  Moreover, the sum of the columns of the inverse is a column vector with strictly positive entries.

\noindent
\linebreak
\textbf{Lemma 2.3.  }  Suppose $L$ is any square matrix with the following properties:  the diagonal entries are positive, the other entries are nonpositive, and in each row the diagonal entry is at least as large as the sum of the absolute values of the off diagonal entries.  In other words, we have
\begin{itemize}
\item{$L_{ii} > 0$,}
\item{$L_{ij} \leq 0$ for $i \neq j$,}
\item{$L_{ii} \geq -\Sigma_{i \neq j} L_{ij}$.}
\end{itemize}
\noindent
In other words, if we sum across any row, we obtain a nonnegative number.  We also require that for each square submatrix (including $L$ itself) obtained by deleting some columns and their respective rows, that one such sum is strictly positive.  Then $L$ is invertible, and its inverse has only nonnegative entries.

This last condition is important in the proof.  It is trivially true for proper submatrices of a Laplacian matrix, provided that the corresponding graph is connected.  For, if such a submatrix had all rows whose entries summed to zero, then no vertex corresponding to one of those rows would be connected to any vertex corresponding to one of the remaining rows.

\noindent
\linebreak
\textit{Proof.  }  The proof is by induction on the size $n$ of $L$.  If $n=1$, the statement is trivial.  In case $n=2$, we have 
$$L = \left( \begin{array}{cc}
a & b \\
c & d \\
\end{array} \right).$$  

\noindent
If $L$ has nonzero determinant, then
$$L^{-1} = \frac{1}{ad-bc} \left( \begin{array}{cc}
d & -b \\
-c & a \\
\end{array} \right) .$$  So we get the result provided that $ad - bc > 0$.  This condition follows from the fact that either $a > b$ or $d >c$.

Suppose that the result is true for matrices of this type with size $n-1$, and let $L$ be such a matrix of size $n$.  We claim that when we clear entries in one of the columns, the resulting $(n-1) \times (n-1)$ matrix obtained by deleting that column and its corresponding row is of the same type as $L$.

Without loss of generality, suppose the first row has entries summing to a strictly positive number (this is true of at least one row).  If there are no entries to clear in the first column, then there is nothing to show.  Otherwise consider one of the rows with a nonzero entry in the first column.  Let this row be $$\left [ \begin{array}{cccc}
a_{j1} & a_{j2} & \cdots & a_{jn} \\
\end{array} \right].$$  Using the same notation for the first row, clearing the leading term from the $j$th row replaces $$a_{jk} \mapsto -\frac{a_{11}}{a_{j1}} a_{jk} + a_{1k}.$$
\noindent
Notice that the $j$th term remains positive.  This amounts to observing that $a_{11}a_{jj} > a_{j1}a_{1j}.$  In fact, we have $$a_{11}a_{jj} \geq (a_{12} + \cdots + a_{1n})(a_{j1} + \cdots + \overline{a_{jj}} + \cdots + a_{jn}),$$ by hypothesis on the rows of $L$, the inequality being strict by hypothesis on the first row.  The off-diagonal entries clearly remain nonnegative.  We have to verify that the sum of the entries in this row is also nonnegative.  But this sum is $$-\frac{a_{11}}{a_{j1}}(a_{j2} + \cdots + a_{jn}) + (a_{12} + \cdots + a_{1n}).$$  Since $$a_{11} > -( a_{12} + \cdots + a_{1n})$$ and $$a_{j2} + \cdots + a_{jn} \geq -a_{j1},$$ the sum of the entries in the new row is positive.

Inductively, the new matrix in the lower right hand corner is invertible, and the inverse has only nonnegative entries.  Since $a_{11}$ is positive, the original matrix is also invertible.  (Row reduction corresponds to multiplication by an invertible matrix.  Such a matrix has nonzero determinant.  After these steps, we can expand along the first column to see that the original matrix also has nonzero determinant.)

We have to show that the inverse of $L$ has only nonnegative entries.  Let us look at what happens to $I_n$ in the augmented matrix $[L|I_n]$ as we row reduce $L$ to the identity.  We start with nonnegative entries only.  Clearing entries in the first column amounted to multiplying a row by a positive number $a_{11}$ and then adding one row to another.  This preserves the property of having only nonnegative entries.  Then, inductively, we row reduce an $(n-1) \times (n-1)$ matrix of the same type to the identity.  By the inductive setup of the proof, all of the steps are the same as in the first step, except possibly we don't start at the first row.  So again, we preserve having nonnegative entries.  Finally, we need to clear entries in the first row and divide by $a_{11}$.  But the off-diagonal entries in the first row are 0 or negative, so we can add positive multiples of the identity columns to clear the nonzero entries.  

More directly, let $B$ be the matrix such that $BL$ has zero entries in the first column except for $a_{11}$.  Denote by $L'$ the $n \times n$ matrix obtained after clearing entries in the first column of $L$. The $n \times n$ matrix $$E = \left( \begin{array}{cccc}
1 & 0 & \cdots & 0 \\
0 & * & * & * \\
\cdots & \cdots & \cdots & \cdots \\
0 & * & * & * \\
\end{array} \right)$$
where the *'s represent the inverse to $L'$, is such that $$EBL = \left( \begin{array}{cccc}
a_{11} & a_{12} & \cdots & a_{1n} \\
0 & * & * & * \\
\cdots & \cdots & \cdots & \cdots \\
0 & * & * & * \\
\end{array} \right)$$
where now the *'s represent the identity matrix.  Finally, if we denote by $A$ the matrix corresponding to the addition of multiples of the bottom rows to the top row of this matrix and then dividing by $a_{11}$, then $AEBL = I_n$.  Since each of $B, E, A$ has nonnegative entries, so does $L^{-1}$.  QED.

\noindent
\linebreak
This invertibility property is directly related to existence and uniqueness of solutions to the problem $$Lu = \chi_{\{u>0\}}.$$
As a consequence, we can determine the reduced echelon form for any Laplacian matrix $L$.  Adding all rows to the bottom row produces a row of zeros.  By the lemma, the remaining $(n-1) \times (n-1)$ matrix in the upper left corner is invertible.  By using an appropriate invertible matrix, then, we can bring $L$ to the form where the bottom row is zero and there is an identity matrix in the upper left corner.  But the vector $u = (1, 1, \ldots, 1)^T$ is such that $Lu = 0$, and this is still true after multiplying $L$ by an invertible matrix $E$.  In other words, $ELu=0$ forces the entries in the last column to be -1, so the reduced row echelon form of $L$ is $$\left( \begin{array}{ccccc}
1 & 0 & \cdots & 0 & -1 \\
0 & 1 & \cdots & 0 & -1 \\
\cdots & \cdots & \cdots & \cdots & \cdots \\
0 & 0 & \cdots & 1 & -1 \\
0 & 0 & 0 & 0 & 0 \\
\end{array} \right).$$
Of course, having found that the rank of $L$ is $n-1$, it follows that 0 is an eigenvalue.  It turns out that 0 is an eigenvalue of multiplicity 1.  Given a Laplacian matrix $L$, its characteristic polynomial can be written as $-t \det L'$, where $L'$ has a first row of all 1's, and the other rows are as in $L - tI$.  Evaluating at 0, we get a matrix with first row all 1's and other rows unaffected.  The proof of Lemma 1.1 shows that the $(n-1) \times (n-1)$ matrix in the lower right corner is invertible.  Thus by multiplying $L'$ by an appropriate invertible matrix, we can bring it to the form where the first row is $$(1,0, 0, \cdots, 0),$$ and there is an identity matrix in the lower right corner.  Such a matrix is lower triangular with 1's on the diagonal, hence invertible.  Thus 0 is not an eigenvalue of the factored matrix $L'$, or equivalently, 0 is an eigenvalue of multiplicity 1.  What are the other eigenvalues?  In the next section, we compute eigenvalues for specific matrices, and discuss a method for finding the characteristic polynomial of any sufficiently simple Laplacian matrix by inspection of the corresponding graph.

\section{Eigenvalues of Laplacians}
\noindent
In this section, we determine the eigenvalues of some specific types of Laplacian matrices.  In particular, we find eigenvalues of matrices corresponding to the cycle graph $C_n$ on $n$ vertices, the path graph $P_n$, the complete graph $K_n$, the complete bipartite graph $K_{m,n}$, and the star of a graph.  At the end of the section, we describe a generalized characteristic polynomial whose coefficients can be determined inductively by inspection of the graph.  As such, its use also seems to be limited to Laplacian matrices.

\noindent
\linebreak
\textbf{The Path Graph.  }  For the path graph on $n$ vertices, the Laplacian matrix is $$P_n = \left( \begin{array}{cccccc}
1 & -1 & 0 & \cdots & 0 & 0 \\
-1 & 2 & -1 & \cdots & 0 & 0 \\
0 & -1 & 2 & \cdots & 0 & 0 \\
\cdots & \cdots & \cdots & \cdots & \cdots & \cdots \\
0 & 0 & 0 & \cdots & -1 & 1 \\
\end{array} \right).$$
\noindent
We show that the eigenvalues are $2 - 2 \cos (\pi k/n)$, $k = 0, 1, \ldots, n-1.$  We do this directly by considering the characteristic polynomial $p(t)$ of $P_n$.  First, we observe that $p(t) = (1-t)D_{n-1}(t) - D_{n-2}(t),$ where $$D_n(t) = \det \left( \begin{array}{cccccccc}
2-t & -1 & 0 & 0 & \cdots & 0 & 0 & 0 \\
-1 & 2-t & -1 & 0 & \cdots & 0 & 0 & 0 \\
0 & -1 & 2-t & -1 & \cdots & 0 & 0 & 0 \\
\cdots & \cdots & \cdots & \cdots & \cdots & \cdots & \cdots & \cdots \\
0 & 0 & 0 & 0 & \cdots & 0 & -1 & 1-t \\
\end{array} \right).$$
\noindent
The recurrence relation for $D_n(t)$ is $D_n = (2-t)D_{n-1}- D_{n-2},$ where $D_1 = 1-t$ and $D_0 = 1.$  The closed form solution of this relation is given by $$D_n(t) = aL^n + b\overline{L}^n,$$ where $$2L = 2-t + \sqrt{t^2-4t}$$ and $$2\overline{L} = 2-t - \sqrt{t^2-4t}.$$  Viewing $t^2-4t$ as imaginary (as a first case, but this ends up being true), we see that $L$ is a complex number with magnitude 1, and that $\overline{L}$ is its conjugate.  When we solve for $a$ and $b$, we obtain $a = 1-b$ and $$b = \frac{L - 1 + t}{L - \overline{L}}.$$
Now $p(t) = D_n(t) - D_{n-1}(t).$  Collecting terms, we see that $p(t) = 0$ when $L^{n-1}(L-1)^2 = \overline{L}^{n-1}(\overline{L}-1)^2.$  We can write $L = e^{i\theta}$.  Using the complex definition of sine and cosine, we are looking for values of $\theta$ for which $$\sin ((n+1)\theta) - 2\sin (n \theta) + \sin ((n-1)\theta)= 0.$$  With the angle sum formulas for sine, this reduces to $$\sin (n\theta)(\cos \theta -1)=0.$$  But when $\cos \theta = 1$, $\sin \theta = \sin (n\theta) =0$, so in any case we obtain $$\theta = 0, \pi/n, 2\pi/n, \ldots.$$  So $p(t) = 0$ when $L = e^{i k \pi/n}$, for $k = 0, 1, \ldots, n-1$.  Finally, using the expressions above for $L$ and identifying real and complex parts, we see that $p(t)$ has zeros $2 - 2 \cos (k \pi/n)$ for $k = 0, 1, \ldots, n-1.$  (Thus in fact the eigenvalues are smaller than 4.)

\noindent
\linebreak
\textbf{The Cycle Graph.  }  For the cycle graph, the characteristic polynomial is $p(t) = F_n(t) - F_{n-2}(t) - 2,$ where $$F_n(t) = \det \left( \begin{array}{cccccc}
2-t & -1 & 0 & 0 & \cdots & 0 \\
-1 & 2-t & -1 & 0 & \cdots & 0 \\
\cdots & \cdots & \cdots & \cdots & \cdots & \cdots \\
0 & 0 & 0 & \cdots & -1 & 2-t \\
\end{array} \right).$$
This is the same matrix as in the previous example, except for the $(n,n)$ entry.  Its recurrence relation is therefore the same relation as before, except that $F_1(t) = 2-t$.  Following the same arguments as above, we see that $p(t) = 0$ when $L = e^{i\theta}$ for $\theta$ satisfying $\cos (n\theta) = 1$; in other words, for $\theta = 0, 2\pi/n, \ldots, 2(n-1)\pi/n$.  Then looking at the real part of $L$, we see that $t$ must equal $2-2\cos (2k\pi/n)$ for $k = 0, 1, \ldots, n-1.$

\noindent
\linebreak
\textbf{The Complete Graph.  }  For the complete graph on $n$ vertices, we are seeking zeros of the polynomial $$p(t) = \det \left( \begin{array}{ccccc}
n-1-t & -1 & -1 & \cdots & -1 \\
-1 & n-1-t & -1 & \cdots & -1 \\
\cdots & \cdots & \cdots & \cdots & \cdots \\
-1 & -1 & -1 & \cdots & n-1-t \\
\end{array} \right).$$
\noindent
Adding all rows to the first row produces a row with all entries $-t$, which can be removed so that the first row is all 1's.  Adding this row to all the other rows does not change the determinant, and we get only diagonal entries each equal to $n-t$.  Thus we have $$p(t) = -t(n-t)^{n-1}.$$  The eigenvalues are 0 and $n$ (with multiplicity $n-1$).

\noindent
\linebreak
\textbf{The Complete Bipartite Graph $K_{m,n}$.  }  We are looking for zeros of the determinant $$p(t) = \det \left( \begin{array}{ccccccccc}
m-t & 0 & 0 & \cdots & 0 & -1 & -1 & \cdots & -1 \\
0 & m-t & 0 & \cdots & 0 & -1 & -1 & \cdots & -1 \\
0 & 0 & m-t & \cdots & 0 & -1 & -1 & \cdots & -1 \\
\cdots & \cdots & \cdots & \cdots & \cdots & \cdots & \cdots & \cdots & \cdots \\
-1 & -1 & -1 & \cdots & -1 & n-t & 0 & \cdots & 0 \\
-1& -1 & -1 & \cdots & -1 & 0 & n-t & \cdots & 0 \\
\cdots & \cdots & \cdots & \cdots & \cdots & \cdots & \cdots & \cdots & \cdots \\
-1 & -1 & -1 & \cdots & -1 & 0 & 0 & \cdots & n-t \\
\end{array} \right),$$
\noindent
where the upper left block is $n \times n$ and the lower right block is $m \times m$.  As before, we add all rows to the first row and then factor out the $-t$.  Then add the resulting row of 1's to the bottom $m$ rows.  We find that $$p(t) = -t(m-t)^{n-1} \cdot \det \left( \begin{array}{ccccc}
n+1-t & 1 & 1 & \cdots & 1 \\
1 & n+1-t & 1 & \cdots & 1 \\
\cdots & \cdots & \cdots & \cdots & \cdots \\
1 & 1 & 1 & \cdots & n+1-t \\
\end{array} \right).$$
\noindent
Performing the same trick on this $m \times m$ matrix, we get $$p(t) = -t(m-t)^{n-1}(m + n - t)(n-t)^{m-1}.$$  So the eigenvalues are 0, $m$ (with multiplicity $n-1$), $m+n$, and $n$ (with multiplicity $m-1$).

\noindent
\linebreak
\textbf{The Star of a Graph.  }  This time we consider what happens to the eigenvalues when we add a vertex that is adjacent to every original vertex.  The Laplacian determinant takes the form $$\det \left( \begin{array}{cccc}
n -t & -1 & \cdots & -1 \\
-1 & * & \cdots & * \\
\cdots & \cdots & \cdots & \cdots \\
-1 & * & \cdots & * \\
\end{array} \right),$$
\noindent
where the *'s represent the Laplacian determinant of the original graph, except that 1 is added to each diagonal entry.  Here $n$ is the original size of the graph before adding the vertex.  We expand along the first column.  The first term produced is $(n-t)p(t-1)$, where $p(t)$ is the Laplacian determinant of the original graph.  The (2,1) minor is, after accounting for signs, $$- \det \left( \begin{array}{cccc}
1 & 1 & \cdots & 1 \\
** & * & \cdots & * \\
\cdots & \cdots & \cdots & \cdots \\
** & * & \cdots & * \\
\end{array} \right),$$
\noindent
where now the *'s represent rows 2 through $n$ of the original matrix, again with 1 added to the diagonal entries.  Note that when finding Laplacian determinants, the first step is to add all rows to the first row and then factor out the $-t$.  The matrix above is what we would get if we had evaluated at $t-1$ and performed the factorization of $-t+1$.  In other words, the (2,1) minor is $$-p(t-1)/(-t+1).$$  Now we note that the (3,1) minor is the same as the (2,1) minor.  The only difference is that we would have added all rows to the second row, but then interchanging rows multiplies the determinant by -1, and this is accounted for by the -1 in the (3,1) entry.  Since there are $n$ terms of this kind, we find our characteristic polynomial $$(n-t)p(t-1) - np(t-1)/(-t+1) = -t(n+1 - t)p(t-1)/(-t+1).$$  Thus we obtain a new eigenvalue $n+1$, and the nonzero eigenvalues increase by 1.

\noindent
\linebreak
\textbf{The Generalized Characteristic Polynomial.  }  The characteristic polynomials for the cycle graph and path graph do not seem, at first, to convey any information about the graphs themselves.  However, there is a polynomial in $n$ variables that does seem to convey useful information about the graph, and from which one can derive the characteristic polynomial in an ordered sequence of steps.  Let $L$ be the Laplacian matrix of a graph with $n$ vertices (of the usual type assumed in this paper).  We define the \textbf{generalized characteristic polynomial} to be the polynomial in $n$ variables $$p(x_1, x_2, \ldots, x_n) = \det \left( L - \left( \begin{array}{ccccc}
x_1 & 0 & 0 & \cdots & 0 \\
0 & x_2 & 0 & \cdots & 0 \\
\cdots & \cdots & \cdots & \cdots & \cdots \\
0 & 0 & 0 & \cdots & x_n \\
\end{array} \right)  \right).$$
\noindent
That is, instead of subtracting $tI$ from $L$, we subtract a diagonal matrix with possibly different entries along the diagonal.  

\noindent
\linebreak
To see why this polynomial is preferable to the usual characteristic polynomial, consider the path graph on $n$ vertices.  We can obtain the usual characteristic polynomial by substituting $x_1 = x_2 = \cdots = x_n = t$ after finding the generalized characteristic polynomial.  We ask what the coefficient of a general term is in this polynomial, say the term $x_{j1}x_{j2}\cdots x_{jl}.$  The determinant in question is $$\det \left( \begin{array}{cccccc}
1-x_1 & -1 & 0 & 0 & \cdots & 0 \\
-1 & 2-x_2 & -1 & 0 & \cdots & 0 \\
\cdots & \cdots & \cdots & \cdots & \cdots & \cdots \\
0 & 0 & 0 & \cdots & -1 & 1 - x_n \\
\end{array} \right).$$
\noindent
To find the term in question, we set $x_k = 0$ for $k \neq j1, j2, \ldots , jl.$  When expanding along the row corresponding to $j1$, the only term that contributes to the desired coefficient is the diagonal entry, and likewise when we expand along rows involving other values of $j$ in subsequent steps.  Notice that deleting a row and column via a diagonal entry creates a block diagonal matrix, at least in the case of the path graph.  We are therefore looking for a product of determinants of two types:  $$\det \left( \begin{array}{ccccccc}
1 & -1 & 0 & 0 & \cdots & 0 & 0 \\
-1 & 2 & -1 & 0 & \cdots & 0 & 0 \\
\cdots & \cdots & \cdots & \cdots & \cdots & \cdots & \cdots \\
0 & 0 & 0 & 0 & \cdots & -1 & 2 \\
\end{array} \right),$$
and 
$$\det \left( \begin{array}{ccccccc}
2 & -1 & 0 & 0 & \cdots & 0 & 0 \\
-1 & 2 & -1 & 0 & \cdots & 0 & 0 \\
\cdots & \cdots & \cdots & \cdots & \cdots & \cdots & \cdots \\
0 & 0 & 0 & 0 & \cdots & -1 & 2 \\
\end{array} \right),$$
the only difference being in the (1,1) entry.  (There is a third type to be computed, with a 2 in the (1,1) entry and a 1 in the $(n,n)$ entry, but this is the same determinant as the first.)  It is straightforward to see that the second determinant is $n+1$, where $n$ is the size of the matrix, and that the first determinant is therefore (recursively) $n - (n-1) = 1.$  It follows that the coefficient of $x_{j1}x_{j2}\cdots x_{jl}$ is $(-1)^l$ times the product of the differences between successive $x_j$ terms.  Likewise, for the cycle graph, the coefficients are products of jumps between successive terms.

\noindent
\linebreak
As an example, we compute the characteristic polynomial for the path graph on 4 vertices.  There is only one term with total degree 4, namely $x_1x_2x_3x_4$, which has coefficient 1.  So the $t^4$ term is just $t^4$.  For the $t^3$ term, there are 4 terms with total degree 3, namely $x_1x_2x_3, x_1x_2x_4, x_1x_3x_4,$ and $x_2x_3x_4.$  The coefficients are $-1, -2, -2,$ and $-1$, respectively.  So the $t^3$ term is $-6t^3$.  The $t^2$ term comes from the degree 2 monomials $x_1x_2, x_1x_3, x_1x_4, x_2x_3, x_2x_4,$ and $x_3x_4$, with respective coefficients 1, 2, 3, 1, 2, and 1, so the $t^2$ term is $10t^2$.  The $t$ term comes from degree 1 monomials, each with coefficient 1, so we get the term $-4t$.  Finally, the degree 0 term is 0 because $p(0,0,\ldots, 0) = 0$ for any Laplacian matrix.  Thus the characteristic polynomial is $$p(t) = t^4-6t^3+ 10t^2 -4t.$$

\noindent
\linebreak
It is not clear, looking at this polynomial or similar polynomials from path graphs with more vertices, how the coefficients are related to the structure of the graph.  The generalized polynomial, admittedly, is easy to compute in the case of a path graph or a cycle graph, but for arbitrary graphs it is difficult to reconcile the coefficients with distances between vertices.  We do have the following result, however.

\noindent
\linebreak
\textbf{Proposition 3.1  }  For an arbitrary graph, label the vertices so that $v_1$ is adjacent to $v_2, v_3, \ldots, v_j$ and no other vertices.  Then the $x_1$-term in $p(x_1, x_2, \ldots, x_n)$ is $-x_1q(x_2-1, x_3-1, \ldots, x_j - 1, x_{j+1}, \ldots, x_n)$, where $q(x_2, x_3, \ldots, x_n)$ is the generalized polynomial for the graph obtained by deleting vertex $v_1$.

\noindent
\linebreak
\textit{Proof.  }  Suppose that the generalized polynomial is given by $$\det \left( \begin{array}{ccccc}
a_{11} -x_1 & a_{12} & a_{13} & \cdots & a_{1n} \\
a_{21} & a_{22}-x_2 & a_{23} & \cdots & a_{2n} \\
\cdots & \cdots & \cdots & \cdots & \cdots \\
a_{n1} & a_{n2} & a_{n3} & \cdots & a_{nn} - x_n \\
\end{array} \right).$$
Expanding along the first row, the only relevant term is the (1,1) minor, and then the term that multiplies $a_{11}$ is the same as the term that multiplies $x_1$.  So the desired term is $-x_1$ times the determinant $$\det \left( \begin{array}{cccc}
a_{22} - x_2 & a_{23} & \cdots & a_{2n} \\
a_{32} & a_{33}- x_3 & \cdots & a_{3n} \\
\cdots & \cdots & \cdots & \cdots \\
a_{n1} & a_{n2} & \cdots & a_{nn} - x_n \\
\end{array} \right).$$
This is not necessarily a Laplacian determinant.  However, the Laplacian determinant for the graph obtained by removing vertex $v_1$ is $$\det \left( \begin{array}{cccc}
a_{22} -1 - x_2 & a_{23} & \cdots & a_{2n} \\
a_{32} & a_{33} -1 - x_3 & \cdots & a_{3n} \\
\cdots & \cdots & \cdots & \cdots \\
a_{n1} & a_{n2} & \cdots & a_{nn}-x_n \\
\end{array} \right),$$
where the first $j-1$ rows have been altered.  The proposition follows from comparing the last two matrices.  QED

\noindent
\linebreak
This result might be useful for building up graphs inductively, adding one new vertex connected to exactly one old vertex, provided that we know the generalized polynomial for the old graph.  One kind of graph where this is useful is a rooted tree.  A 1-level rooted tree is just a star graph (a complete bipartite graph $K_{m,1}$).  An arbitrary rooted tree is obtained by adding one vertex at a time, connected to exactly one older vertex.  As it happens, this method of finding the characteristic polynomial for the case of a tree seems to be extremely inefficient.  However, the method does bring out the fact that the coefficients of the characteristic polynomial are related to the structure of the graph.

\noindent
\linebreak
\textbf{Example 3.2.  }  As a simple example, we give an expression for the characteristic polynomial of the graph obtained by adding a new edge at the second vertex of a path graph on $n$ vertices.  (This is sometimes called the graph $D_n$, whereas the path graph is also called $A_n$.)

\noindent
\linebreak
\textit{Solution.  }  Looking at the matrices in Proposition 3.1, we see that adding a vertex $v_{n+1}$ adjacent to vertex $v_2$ creates the generalized characteristic polynomial $$(1-x_{n+1})p_n(x_1, x_2-1, x_3, \ldots , x_n) - (1-x_1)p_{n-2}(x_3-1, x_4, \ldots, x_n).$$  Here we use $p_n$ to denote the generalized characteristic polynomial for the path on $n$ vertices, and we have renamed the variables so they correspond to labeling of the original path graph.  The polynomials $p_n$ can be determined (somewhat painstakingly) by the method described above.  Finally, evaluating this expression at $x_1 = x_2 = \cdots = x_{n+1} = t$ gives the usual characteristic polynomial.  

\noindent
\linebreak
In future research, we intend to address the problem of adding an edge to an existing (non-complete) graph on $n$ vertices.  Together with a full understanding of how the characteristic polynomial changes upon adding a vertex and a single edge, it may be possible to efficiently determine the generalized characteristic polynomial by inspection of the graph alone.

In the final section of the paper, we return to the motivating problem, the obstacle problem on graphs.  The computation of eigenvalues is not directly related to this problem.

\section{The Obstacle Problem on Graphs}  
\noindent
In this section we prove Theorems 2.1 and 2.2 and give some examples.

\noindent
\linebreak
\textbf{Theorem 4.1.  }  For any Laplacian matrix $L$, there is a unique solution to the obstacle problem $$Lu = \chi_{\{u>0\}}.$$   That is, for each assignment of the free boundary $u_i = 0$ for $i \in I$, where $I \subset \{1, 2, \ldots, n\}$ is nonempty, there is a unique positive value for $u_i$, $i \notin I$, satisfying the nonzero equations in the system $Lu = \chi_{\{u>0\}}.$

\noindent
\linebreak
\textit{Proof.  }  This is essentially Lemma 2.3 again.  Since $I$ is nonempty, we are considering a proper subset of the set of linear expressions defined by $L$.  By setting some of the $u_i = 0$, we omit the corresponding equations from the given system.  The resulting system of equations, each with right hand side 1, is obtained by considering a proper submatrix of the Laplacian, which we have established as invertible.  If this matrix is $L'$ and the nonzero coordinates of $u$ are written as the vector $u'$, then we obtain $u' = L'^{-1}\overline{1},$ where $\overline{1}$ represents the vector consisting entirely of 1's.  Thus there exists a unique solution, and we have to establish that it is strictly positive.  As the proof of Lemma 2.3 shows, the formation of the inverse matrix $L'^{-1}$ involves changing the identity matrix by multiplying rows by positive constants and adding rows together.  In any case, the sum of the columns remains strictly positive under each operation.  This sum is precisely $L'^{-1}\overline{1}.$  QED.

\noindent
\linebreak
\textbf{Theorem 4.2  }  For an arbitrary system $Lu + b = \chi_{\{u>0\}}$, where $b$ is a random vector (with the properties given in the Introduction), the nonzero entries of $b$ and the nonzero values $u_i$ are uniquely determined and positive.

\noindent
\linebreak
\textit{Proof.  }  In fact, the nonzero $u_i$ are determined as in Theorem 4.1, hence are uniquely determined and positive.  In this case, we do consider the other equations corresponding to vertices where $u_i = 0$.  Since the values of the $u_i$ in these equations are already determined and positive, the values $b_i$ are also uniquely determined and positive.  QED.

\noindent
\linebreak
\textbf{The Complete Graph.  }  Without loss of generality, assume that $u_1 = u_2 = \cdots = u_j = 0.$  The first $j$ equations in the system have the same left hand side:  $$-u_{j+1} - u_{j+2} - \cdots - u_n$$ and the other equations are of the type $$(n-1)u_{j+1} - u_{j+2} - \cdots - u_n = 1.$$  Subtracting any two of the second type of equations yields $u_{j+1} = u_{j+2} = \cdots = u_n$.  This common value is determined by $((n-1)-(n-j-1)))u_{j+1} = 1,$ so that $$u_{j+1} = \cdots = u_n = 1/j.$$  In case we introduce an error term, then its nonzero coordinates are all equal to $(n-j)/j$. 

\noindent
\linebreak
Are there other ways of interpreting the usual obstacle problem $Lu= \chi_{\{u>0\}}$?  We have on the one hand ignored auxiliary equations, and on the other modified those equations so that they are true by definition.  It turns out that there is a middle of the road approach.  Consider the system $$Lu + b = \chi_{\{u>0\}},$$ where now $b$ is a constant vector.  

Adding a constant accounts for the fundamental inconsistencies in the original system.  In fact, adding all the equations together specifies what the value of $b$ is for each assignment of zeros:  $b = (n-j)/n,$ where $j$ is the number of zero terms $u_i$.   In the present case, we are not going so far as to say that extra equations are automatically true.  The system has been modified, but now all equations are significant.  The trouble is that solutions do not always exist, though they do exist uniquely and are positive, when they exist.  Moreover, we obtain a scalar multiple by $1-b$ of the solution $u$ found by either of the above methods, provided that we get a solution.

Finding assignments where solutions exist depends on the symmetries of the graph being considered.  For the complete graph, there is always a (unique, positive) solution.  For the cycle graph, for any sufficiently large number of zeros, their assignments to the graph must be symmetrically placed, and likewise for the path graph.  As an example of how this approach works, we consider the complete bipartite graph.

\noindent
\linebreak
\textbf{The Complete Bipartite Graph.  }
Label the vertices $v_1, v_2, \ldots, v_m$ on one side so that $v_1 = v_2 = \cdots = v_r = 0$, and $w_1, w_2, \ldots, w_n$ on the other side so that $w_1 = w_2 = \cdots = w_s = 0$.  The system of equations reduces to $$\begin{array}{c}
-w_{s+1} - w_{s+2} - \cdots - w_n = -b\\
-v_{r+1} - v_{r+2} - \cdots - v_n = -b\\
nv_{r+1} - w_{s+1} - \cdots - w_n = 1-b \\
mw_{s+1} - v_{r+1} - \cdots - v_m = 1-b, \\
\end{array}$$
where there are other equations of the last two types corresponding to the positive values of $v$ and $w$.  We obtain $v_{r+1} = v_{r+2} = \cdots = v_m$ and $w_{s+1} = \cdots = w_n,$ by subtracting like equations of the last two types.  The first two equations then yield the common values $v_m = b/(n-r)$ and $w_n = b/(n-s)$, so the solution is uniquely determined and positive, if it exists.    

In order to have a consistent system, the last two equations have to make sense.  We require that $nv_m - (n-s)w_n = 1-b$ and that $mw_n - (m-r)v_m = 1-b.$  In other words, $(m+n-s)w_n = (n+m-r)v_m$.  With the above calculation, this is true if and only if $1 + m/(n-r) = 1 + m/(n-s),$ or that $r=s$.  

Adding the equations together, we find that $b = (m+n-r-s)/(m+n)$.  The last two equations reduce to $sv_m = 1-b.$  Since $b = 1 - 2s/(m+n)$, it must also be true that $$2/(m+n)= b/(n-s) = 1/(n-s) -2s/((m+n)(n-s)).$$ or $$2(n-s) = m+n -2s,$$ or that $m=n$.   Thus the system has a solution only when $m = n$ and $r=s$.

\noindent
\linebreak
Future research will involve other functions' being added to the $Lu$ side besides a constant or a random vector.  However, it is worth noting that any function of $\chi_{\{u>0\}}$ is necessarily of the type $a\chi_{\{u>0\}} + b$ for constants $a$ and $b$.  So changing the left hand side must involve new functions of $u$.  For example, it might be interesting to find solutions of the problem $Lu + f(u) = \chi_{\{u>0\}}$ for some function $f$.  It seems inevitable, though, that solutions of the obstacle problem on graphs should be determined only when certain consistencies hold with respect to equations corresponding to vertices at which $u_i = 0$.

\section{References}  

\noindent
\linebreak
 1.  [BT]  Ivan Blank and Kubrom Teka, ``The Caffarelli Alternative in Measure for the Nondivergence Form Elliptic Obstacle Problem with Principal Coefficients in VMO."  \textit{Communications in Partial Differential Equations}, Vol. 39, Issue 2, 2014.
 
\noindent
\linebreak
2. [M]  Bojan Mohar, ``Some Applications of Laplace Eigenvalues of Graphs,"  University of Ljubljana.  Also appeared in ``Graph Symmetry:  Algebraic Methods and Applications," Eds. G. Hahn and G. Sabidussi, NATA ASI Ser. C 497, Kluwer, 1997, pp. 225-275.
 
\noindent
\linebreak
3.  [PSU]  Arshak Petrosyan, Henrik Shahgholian, and Nina Uraltseva, \textit{Regularity of Free Boundaries in Obstacle-Type Problems}, Graduate Studies in Mathematics, Vol. 136,  American Mathematical Society, 2012.

\noindent
\linebreak
4.  [S]  Daniel Spielman, Yale University, Lecture Notes on Spectral Graph Theory, 2009.

\end{document}